\documentclass[11pt,reqno]{amsart} 

\usepackage{amssymb,latexsym}
\usepackage{cite} 

\usepackage{fontenc}
\usepackage{times}

\usepackage{amsfonts,amssymb,epsfig}
\usepackage{amscd}
\usepackage{color}
 \usepackage{hyperref}

\usepackage{textcomp}
\usepackage{mathrsfs}

\usepackage[height=200mm,width=140mm]{geometry} 

\theoremstyle{plain}

\theoremstyle{remark}

\theoremstyle{example}

\numberwithin{equation}{section} 

\theoremstyle{definition}
\newtheorem{definition}{Definition}[section]

\newtheorem{example}[definition]{Example}
\newtheorem{remark}[definition]{Remark}

\newtheorem{theorem}[definition]{Theorem}
\newtheorem{proposition}[definition]{Proposition}
\newtheorem{corollary}[definition]{Corollary}
\newtheorem{lemma}[definition]{Lemma}

\newenvironment{prof1}[1][Proof of the Proposition \ref{Propo2}]{\textbf{#1.} }{\ \rule{0.5em}{0.5em}}
\newenvironment{prof2}[1][Proof of the Theorem \ref{Theorem12}]{\textbf{#1.} }{\ \rule{0.5em}{0.5em}}
\newenvironment{proof1}[1][Proof]{\textbf{#1.} }{\ \rule{0.5em}{0.5em}}
\begin{document}
\title[Finslerian convolution metrics and their special classes]{Finslerian convolution metrics and their special classes} 

\author{Gilbert Nibaruta and Esp\'erance Niyonkuru}

\address{Ecole Normale Sup\'erieure\\ 
\noindent
P.~O. Box 6983 \\ 
 Bujumbura\\ Burundi}

 \email{nibarutag@gmail.com}
 
 \address{Institut de P\'edagogie Appliqu\'ee\\ 
 	\noindent
 	D\'epartement de Math\'ematiques\\
 	Bujumbura\\ Burundi}
 
 \email{esperanceni59@gmail.com}
\begin{abstract}
Here, it is introduced a concept of convolution metric in Finslerian Geometry. This convolution metric is a kind of function obtained by a given mathematical operation between two Finslerian metrics. 
Some basic properties of the Finslerian convolution metrics are studied. Then it is characterized Finslerian convolution metrics which are of type Riemannian, Minkowskian as well as Randers. Furthermore, some examples of the Finslerian convolutions are given.
\end{abstract} 
\subjclass[2010]{primary 53B40, secondary 58B20}

\keywords{Finslerian metric; Riemannian metric; Randers metric; convolution metric.}

\maketitle

\section{Introduction}\label{Section1}
Convolution metrics play an important role in Differential Geometry and have applications that include General Relativity since they are natural extensions of warped product metrics. As a concrete example of application, many basic solutions of the Einstein field equations on Riemannian manifolds have solutions when the Riemannian metrics are warped products \cite{Kuhnel} at most.
  
Nowadays, a Finslerian metric is of considerable interest due to its rich structure
  including mainly Riemann, Randers, Minkowski and Berwald type metrics. For the details, the study of each of these special Finslerian metrics and some
important examples can be found in \cite{bao,Nibaruta9} and in \cite{Bin}.
  
   In Riemannian Geometry, the notion of convolution metrics was introduced by Bang-Yen Chen to establish and characterize the Euclidean version of Segre embedding \cite{bang2002}. The following Proposition shows that the notion of convolution metrics arises very naturally from tensor product immersions.
   \begin{proposition}\cite{bang2016}
    Let $M_1$ and $M_2$ be two $C^{\infty}$ manifolds of dimensions $n_1$ and $n_2$ respectively. Consider $g_1$ and $g_2$ two Riemannian metrics on $M_1$ and $M_2$ respectively. If $\varphi_1:(M_1,g_1)\longrightarrow \mathbb{R}^{n_1}\backslash \{0\}$ and $\varphi_2:(M_2,g_2)\longrightarrow \mathbb{R}^{n_2}\backslash \{0\}$ are isometric  immersions then the mapping 
    \begin{eqnarray}\varphi=\varphi_1\otimes\varphi_2:M_1\times M_2\longrightarrow \mathbb{R}^{n_1}\otimes\mathbb{R}^{n_2}=\mathbb{R}^{n_1n_2}\nonumber\end{eqnarray} 
    defined by $\varphi(x_1,x_2)=\varphi_1(x_1)\varphi_2(x_2)$ gives rise to a convolution manifold $M_1\bigstar M_2$ endowed with the convolution metric 
    \begin{eqnarray}_{\rho_2}g_1\star_{\rho_1}g_2=\rho_2^2g_1+\rho_1^2g_2+2\rho_1\rho_2d\rho_1\otimes d\rho_2\nonumber\end{eqnarray}
    whenever $_{\rho_2}g_1\star_{\rho_1}g_2$ is non-degenerate where $\rho_1=\| \varphi_1\|$ and $\rho_2=\| \varphi_2\|$.
   \end{proposition}
   
   In this work, we construct an extension of the Riemannian convolution metric on a product manifold (See the Definition \ref{Defi03}). More precisely, we define a new notion of convolution metric by taking account of, not only the Riemannian metrics but also the Minkowski, the Randers and the Berwald metrics.
  We have the following: 
  \begin{proposition}\label{Propo2}
 Let $(M_1,F_1)$ and $(M_2,F_2)$ be two Finslerian manifolds. Consider $f_k:M_k\longrightarrow (0,\infty)$, for $k\in \{1,2\}$, some smooth functions. 
If $_{f_2}F_1\star_{f_1}F_2$ is a convolution Finslerian metric on a manifold $M_1\times M_2$ then one of the following assertions holds
\begin{itemize}
\item[(1)]the gradient of one of the functions $f_k$ is orthogonal to any tangent vector on $M_1\times M_2$.
\item[(2)] one of the $f_k$ is a nonzero constant. In this last case, when one of the $f_k\equiv 1$, the convolution Finslerian metric  $_{f_2}F_1\star_{f_1}F_2$ reduces to a Finslerian warped product metric.
\end{itemize} 
\end{proposition}
\noindent
   Main results regard the characterization of convolution metrics of type Riemannian, Minkowskian as well as Randers among Finslerian ones. For example, we prove the following.
   \begin{theorem}\label{Theorem12}
	Let $(M_1,F_1)$ and $(M_2,F_2)$ be two Finslerian  manifolds of Randers type that is, for $k\in\{1,2\}$, $F_k:=\alpha_k+\beta_k$ where
	$\alpha_k$ is a Riemannian metric 
	and
	$\beta_k$ is a $1$-form on $M_k$ with $||\beta_k||_{\alpha_k}<1$. If one of the $f_k:M_k\longrightarrow (0,\infty)$ is a nonzero constant function, then the Finslerian convolution metric $_{f_2}F_1\star_{f_1}F_2$ on $M_1\times M_2$ is a Randers metric if and only if 
	\begin{eqnarray}\label{34}
	\frac{\alpha_1}{\alpha_2}=\frac{\beta_1}{\beta_2}.
	\end{eqnarray}
\end{theorem}
   \noindent
   One of the examples of special Finslerian convolution is given by
   \begin{example}\label{example11}
 Consider $M=\{x=(x^1,x^2,x^3,x^4)\in\mathbb{R}^4:x_1>0,x_3>0\}$ and $T_xM=\{(y^1,y^2,y^3,y^4)\in\mathbb{R}^4:y_1>0,y_3>0\}$. For $\lambda\in [2,4]$ and for a positive integer $k$, define $F:TM\longrightarrow (0,\infty)$ by 
 \begin{eqnarray}\label{dhggdgdgh}
    F(x^1,x^2,x^3,x^4,y^1,y^2,y^3,y^4)
    &=&\left\{x_3^2\left(y_1^4+\lambda y_1^2y_2^2+y_2^4\right)^{\frac{1}{2}}
    +8x_1^3x_3^3y_1y_3\right.\nonumber\\
    &&+\left. x_3^2\left[y_3^2+y_4^2+\lambda \left(y_3^{2k}+y_4^{2k}\right)^{\frac{1}{k}}\right]\right\}^{\frac{1}{2}}.
\end{eqnarray}
Then the function $F$ defined in (\ref{dhggdgdgh}) is a Finslerian convolution metric.
\end{example}

 The outline of this document is organized as follows. 
  In Section \ref{Section2},  
  we establish some basic notions on Finslerian convolution manifolds.
  The Section \ref{Section3} is devoted to the characterizations of the convolution of Riemannian, Minkowskian and Randers metrics among the Finslerian convolution metrics. In the Section \ref{Section4} we provide examples of Finslerian convolution metrics.\\
  
  \noindent
\textbf{Notation.}
 In the following, $M$ is an $n-$dimensional product of the $C^2$-manifolds $M_1$ and $M_2$ of dimensions $n_1$ and $n_2$ respectively. We denote by $T_xM$ the tangent space at 
   $x=(x_1,x_2)\in M$ and by $TM:=\bigcup_{x\in M}T_xM$ the tangent bundle of $M$.
   Set $\mathring{T}M=\{(x_1,x_2,y_1,y_2)\in TM: y_1\neq 0, y_2\neq 0 \text{ where } 0 \text{ is a zero vector } \in T_{x_k}M_k \text{ for } k\in\{1,2\}\}$ and $\pi:TM\longrightarrow M:\pi(x,y)\longmapsto x$ the natural projection. 
   Let $(x^1,...,x^n)=(x_1^1,...,x_1^{n_1},x_2^{n_1+1},...,x_2^{n})$ be a local coordinate 
   on an open subset $U$ of $M$ and $(x^1,...,x^n,y^1,...,y^n)$ be the local coordinate 
   on $\pi^{-1}(U)\subset TM$. 
   The local coordinate system $(x^i)_{i=1,...,n}$ produces the coordinate bases  
   $\{\frac{\partial}{\partial x^i}\}_{i=1,...,n}$ and $\{dx^i\}_{i=1,...,n}$ respectively, for 
   $TM$ and cotangent bundle $T^*M$.\\
  
  \noindent
\textbf{Convention.}
   We use Einstein summation convention: repeated
   upper and lower indices will automatically be summed, from $1$ to $n$, unless otherwise will be noted. For example 
   $g_{ij}v^iv^j=\sum_{i,j=1}^n g_{ij}v^iv^j.$
  
 \section{Finslerian convolution manifold}\label{Section2} 

   \begin{definition}\label{defi1}Let $M$ be an $n$-dimensional manifold. A function $F:TM\longrightarrow [0,\infty)$ is a \textit{Finslerian metric} 
   	on $M$ if :
   	\begin{itemize}
   		\item [(i)] $F$ is $C^{\infty}$ on the entire slit tangent bundle $\mathring{T}M$,
   		\item [(ii)] $F$ is positively $1$-homogeneous on the fibers of $TM$, that is \\
   		$\forall c>0,~F(x,cy)=cF(x,y),$
   		\item [(iii)] the Hessian matrix $(g_{ij}(x,y))_{1\leq i,j\leq n}$ with elements
   		\begin{eqnarray}\label{1}
   		g_{ij}(x,y):=\frac{1}{2}\frac{\partial^2F^2(x,y)}{\partial y^i\partial y^j}
   		\end{eqnarray}
   		is positive-definite at every point $(x,y)$ of $\mathring{T}M$.
   	\end{itemize}
   \end{definition}
\begin{remark}
	The Hessian matrix $g$ whose elements are defined in (\ref{defi1}) is a natural Riemannian metric on the pulled-back bundle $\pi^*TM$ over the manifold 
	$\mathring{T}M$. $\pi^*TM$ is a vector bundle and its fiber at a point $(x,y)\in\mathring{T}M$ is
\begin{eqnarray}\label{12}
 \pi^*TM |_{(x,y)}:=\{(x,y;v): v\in T_xM\}|_{(x,y)}\cong T_xM.
\end{eqnarray}
\end{remark}
\begin{definition}\label{defi2}
	Let $(M,F)$ be an $n$-dimensional Finslerian manifold and $g$ the Hessian matrix associated with $F$.
	The gradient of a smooth function $u$ on $M$ is, the section of the vector bundle $\pi^*TM$ denoted by  $\nabla u$, 
	 given by
	 \begin{eqnarray}
	 	g_{(x,y)}(\nabla u_{(x,y)},\xi_{(x,y)})=du_{(x,y)}(\xi_{(x,y)})
	 \end{eqnarray}
	 for any $\xi\in\Gamma(\pi^*TM)$ and for every $(x,y)\in \mathring{T}M$. Locally, 
	 \begin{eqnarray}
	 \nabla u_{(x,y)}=g^{ij}(x,y)\frac{\partial u}{\partial x^i}\frac{\partial }{\partial x^j}.
	 \end{eqnarray}
\end{definition}
\begin{remark}
Let $M_1$ and $M_2$ be two $C^{\infty}$ manifolds. 
For every $(x_1,x_2)\in M_1\times M_2$, we have the following properties 
derived from $M_1$ and $M_2$.
\begin{itemize}
	\item [(1)] The projections 
	\begin{eqnarray}
	\sigma_1&:&M_1\times M_2\longrightarrow M_1
	\text{		such that		} \sigma_1(x_1,x_2)=x_1\nonumber\\
	\sigma_2&:&M_1\times M_2\longrightarrow M_1
	\text{		such that		} \sigma_2(x_1,x_2)=x_2\nonumber
	\end{eqnarray}
	are $C^{\infty}$ submersions.
	\item[(2)] $dim(M_1\times M_2)=dim 
	M_1+ dim M_2$.
\end{itemize}
\end{remark}
   \begin{definition}\label{Defi03}
   	Let $(M_1,F_1)$ and $(M_2,F_2)$ be two Finslerian manifolds. Consider $f_k:M_k\longrightarrow (0,\infty)$, for $k\in\{1;2\}$, some smooth functions. A convolution of $(M_1,F_1)$ and $(M_2,F_2)$, via $f_1$ and $f_2$, is the product manifold $M_1\times M_2$ endowed with a function $F_{f_1,f_2}:T(M_1\times M_2)\longrightarrow (0,\infty)$ defined by
   	\begin{eqnarray}\label{4}
   	F_{f_1,f_2}=
   	\sqrt{h_2^2(\sigma_1^*F_1)^2
   		+h_1^2(\sigma_2^*F_2)^2
   		+\frac{1}{2}h_1h_2\frac{\partial (\sigma_1^*F_1)^2}{\partial y^i}\frac{\partial (\sigma_2^*F_2)^2}{\partial y^j}\big(\nabla h_1\big)^i\big(\nabla h_2\big)^j}
   	\end{eqnarray}
   	where $h_k:=f_k\circ\sigma_k$.
   	The function $F_{f_1,f_2}$ is called a Finslerian convolution of $F_1$ and $F_2$, via $f_1$ and $f_2$, and is denoted by $_{f_2}F_1\star_{f_1}F_2$. Given a Finslerian convolution function  $_{f_2}F_1\star_{f_1}F_2$ on $M_1\times M_2$, the pair $(M_1\times M_2,_{f_2}F_1\star_{f_1}F_2)$ is called a Finslerian convolution
    manifold and is denoted by $_{f_2}M_1\bigstar_{f_1}M_2$.
   \end{definition}
\begin{remark}For any $x\equiv (x_1,x_2)\in M\equiv M_1\times M_2$ and $y\equiv (y_1,y_2)\in T_xM$,
	\begin{itemize}
		\item [(1)] $F_{f_1,f_2}(x,y)$ exists for some positive-real functions $f_1$ and $f_2$ on $M_1$ and $M_2$ respectively. 
		\item [(2)]$F_{f_1,f_2}$ is not $C^{\infty}$ on the tangent vectors 
		of the form $(y_1,0)$ nor $(0,y_2)$ at a point $(x_1,x_2)\in M_1\times M_2$ but on $(y_1,y_2)$ for $y_1\neq 0$ and $y_2\neq 0$.
		\item[(3)] $_{f_2}F_1\star_{f_1}F_2$ is positively homogeneous of degree one with respect to $(y_1,y_2)$.
		\item [(4)] from the relation (\ref{4}), we have 
		{\small{\begin{eqnarray}\label{5}
		F^2_{f_1,f_2}(x,y)
		&=&2f_1(x_1)f_2(x_2)F_1(x_1,y_1)F_2(x_2,y_2)\frac{\partial F_1(x_1,y_1)}{\partial y^i}\frac{\partial F_2(x_2,y_2)}{\partial y^j}\big(\nabla f_1(x_1)\big)^i\big(\nabla f_2(x_2)\big)^j
		\nonumber\\
		&&+f^2_2(x_2)F^2_1(x_1,y_1)+f_1(x_1)^2F^2_2(x_2,y_2).
		\end{eqnarray}}}
		Since $(x_1,x_2)$ is arbitrary a point of $M$, the expression (\ref{5}) can be written as 
		\begin{eqnarray}\label{6}
		F^2_{f_1,f_2}=f_2^2F_1^2+f_1^2F_2^2+2f_1f_2F_1F_2\frac{\partial F_1}{\partial y^i}\frac{\partial F_2}{\partial y^j}\big(\nabla f_1\big)^i\big(\nabla f_2\big)^j.
		\end{eqnarray}
	\end{itemize}
\end{remark}
\begin{lemma}\label{lem11}
Let $M_1$ and $M_2$ be two $C^{\infty}$ manifolds of dimensions $n_1$ and $n_2$, respectively, with $n_1+n_2=n$.  If $F=_{f_2}F_1\star_{f_1}F_2$ is a Finslerian convolution function on the manifold $M_1\times M_2$ then the Hessian matrix associated with $F$ is given by 
\begin{eqnarray}\label{3.1}
\Big(g_{ij}(x_1,x_2,y_1,y_2)\Big)=
\left(\begin{array}{cc}f_2^2(x_2)\Big(g_{1\alpha\beta}(x_1,y_1)\Big) & 2f_1(x_1)f_2(x_2)\Big(\frac{\partial f_1(x_1)}{\partial x_1^{\alpha}}\frac{\partial f_2(x_2)}{\partial x_2^{\nu}}\Big)\\
\Big(0\Big)  & f_1^2(x_1)\Big(g_{2\mu\nu}(x_2,y_2)\Big)
\end{array} \right)\nonumber
\end{eqnarray} 
at every point $(x_1,x_2,y_1,y_2)$ of $\mathring{T}M_1\times\mathring{T} M_2$
where $i,j\in [1,...,n], \alpha,\beta\in [1,...,n_1], \mu,\nu\in[n_1+1,...,n]$ and, $g_1$ and $g_2$ are the Hessian matrices of $F_1$ and $F_2$ respectively.
\end{lemma}
\begin{proof1}
The proof follows from a straightforward calculation using the Definition \ref{defi1} and the relation (\ref{6}).
\end{proof1}
\begin{remark}
 In general, the components $g_{ij}$ of the convolution matrix $g$ are not symmetric in all indices.
\end{remark}
\begin{proposition}\label{Proposition 2.2}
Let $M_1$ and $M_2$ be two $C^{\infty}$ manifolds of dimensions $n_1$ and $n_2$, respectively, with $n_1+n_2=n$. If $_{f_2}F_1\star_{f_1}F_2$ is a Finslerian convolution metric on a manifold $M_1\times M_2$ and if $v\neq 0$ is a  tangent vector on $M_1\times M_2$ then the quadratic form $(_{f_2}g_1\star_{f_1}g_2)(v,v)$ associated with $F$ is positive-definite if
\begin{eqnarray}\label{2.7}
 \frac{1}{f_1^2}g_{1\alpha\beta}v^{\alpha}v^{\beta}+\frac{1}{f_2^2}g_{2}v^{\mu}v^{\nu}
 >-\frac{2}{f_1f_2}\frac{\partial f_1}{\partial x_1^{\alpha}}\frac{\partial f_2}{\partial x_2^{\nu}}
 v^{\alpha}v^{\nu}
\end{eqnarray}
where $\alpha,\beta\in [1,...,n_1], \mu,\nu\in[n_1+1,...,n]$ and, $g_1$ and $g_2$ are the Hessian matrices of $F_1$ and $F_2$ respectively.
\end{proposition}
\begin{proof1}
 From the Lemma \ref{lem11}, for any tangent vector $v\neq 0$ on $M_1\times M_2$, we have
 \begin{eqnarray}
  g(v,v)
  &=&\sum_{i,j=1}^n g_{ij}v^iv^j\nonumber\\
  &=&\sum_{\alpha,\beta=1}^{n_1} g_{\alpha\beta}v^{\alpha} v^{\beta}
  +\sum_{\alpha=1}^{n_1} \sum_{\mu=n_1+1}^{n}g_{\alpha\mu}v^{\alpha}v^{\mu}
  +\sum_{\mu=n_1+1}^{n}\sum_{\alpha=1}^{n_1} g_{\mu\alpha}v^{\mu}v^{\alpha}
  +\sum_{\mu,\nu=n_1+1}^{n} g_{\mu\nu}v^{\mu} v^{\nu}.\nonumber
 \end{eqnarray}
 Since $g\equiv _{f_2}g_1\star_{f_1}g_2$, it follows that
 \begin{eqnarray}\label{2.9}
  g_{ij}(v,v)
  &=&\sum_{\alpha,\beta=1}^{n_1} f_2^2 g_{1\alpha\beta}v^{\alpha} v^{\beta}
  +\sum_{\alpha=1}^{n_1} \sum_{\mu=n_1+1}^{n}2f_1f_2
  \frac{\partial f_1}{\partial x_1^{\alpha}}\frac{\partial f_2}{\partial x_2^{\mu}}
  v^{\alpha}v^{\mu}
  +0
  +\sum_{\mu,\nu=n_1+1}^{n} f_1^2g_{\mu\nu}v^{\mu} v^{\nu}\nonumber\\
 \end{eqnarray}
 with $g_{\alpha\mu}:=2f_1f_2
  \frac{\partial f_1}{\partial x_1^{\alpha}}\frac{\partial f_2}{\partial x_2^{\mu}}$.
  
  By the fact that $F_1$ and $F_2$ are Finslerian metrics, the first and the last terms in the right hand-side of (\ref{2.9}) are positives. Hence $g(v,v)>0$ if the equation (\ref{2.7}) holds.
\end{proof1}

When the matrix given in the Lemma \ref{lem11} is positive-definite symmetric on $\mathring{T}(M_1\times M_2)$, it defines a Finslerian metric
on $M_1\times M_2$. In this case, $_{f_2}F_1\star_{f_1}F_2$ is called a convolution Finslerian metric and the Finslerian convolution
manifold $_{f_2}M_1\bigstar_{f_1}M_2$ is called a convolution Finslerian manifold.

\begin{prof1}
The proof is obtained by the Lemma \ref{lem11} together with the Proposition \ref{Proposition 2.2}.
\end{prof1}

\section{Important results}\label{Section3}
\begin{definition}\label{Defi04}
A Riemannian metric $g$ on an $n$-dimensional manifold $M$
is a family $\{g_x\}_{x\in M}$ where $g_x$ is an inner product on the tangent space $T_xM$, 
such that the elements $
g_{ij}(x):=g_x\left(\frac{\partial}{\partial x^i},\frac{\partial}{\partial x^j}\right) 
$
are $C^{\infty}$ in local coordinates. Since $g_x$ is an inner product, the matrix $(g_{ij}(x))$ is positive-definite at every point $x\in M$.
Thus every Riemannian metric $g$ is a Finslerian metric 
which arises in the following manner
\begin{eqnarray}
F(x,y):=\sqrt{g_{ij}(x)y^iy^j}.\label{7}
\end{eqnarray}
In that case, $\frac{1}{2}\frac{\partial^2F^2(x,y)}{\partial y^i\partial y^j}$ is simply $g_{ij}(x)$ which is independant of $y$.
\end{definition}
\begin{theorem}
	A convolution Finslerian metric $_{f_2}F_1\star_{f_1}F_2$ is Riemannian  if and only if $F_1$ as well as $F_2$ are Riemannian metrics. 
\end{theorem}
\begin{proof1}
	Assume that $F=_{f_2}F_1\star_{f_1}F_2$ is a Riemannian metric. By the Definition \ref{Defi04} the components $g_{ij}$ associated with $F$ are independent of $(y_1,y_2)$. Hence, from the Lemma \ref{lem11}, we have $$g_{1\alpha\beta}(x_1,y_1)=g_{1\alpha\beta}(x_1) \text{      and      } g_{2\mu\nu}(x_2,y_2)=g_{2\mu\nu}(x_2).$$
	That is $F_1(x_1,y_1)$ and $F_2(x_2,y_2)$ are Riemannian metrics. 
	
	The converse is obvious.
\end{proof1}
\begin{definition}
	A Finslerian metric $F$ on a manifold $M$ is locally Minkowski metric if, 
	for every $x\in M$, there exists an neighborhood $V$ of $x$ such that the components $g_{ij}$
	of the Hessian matrix, given in (\ref{1}), satisfy 
	\begin{eqnarray}
	g_{ij}(x,y)=g_{ij}(y).\label{fdkjfjdf}
	\end{eqnarray}
\end{definition}
\begin{theorem}\label{theorem32}
	A convolution Finslerian metric $_{f_2}F_1\star_{f_1}F_2$ is locally Minkowskian  if and only if $F_1$ as well as $F_2$ are locally Minkowskian metrics and the functions $f_k,k\in\{1,2\}$ are both nonzero constants.
\end{theorem}
\begin{proof1}
 The Theorem \ref{theorem32} is obtained by a short calculation of the locally Minkowskian metric from the Lemma \ref{lem11} by using the expression (\ref{fdkjfjdf}).
\end{proof1}

\begin{definition}\label{defi33}
	A Randers metric $F$ on an $n$-dimensional manifold $M$  is a Finsler metric 
	which has the following form 
	\begin{eqnarray}\label{33}
	F:=\alpha+\beta
	\end{eqnarray}
	where, for every $(x,y)\in \mathring{T}M$, 
	$\alpha(x,y):=\sqrt{a_{ij}(x)y^iy^j}$ is a Riemannian metric 
	and
	$\beta(x,y):=b_i(x)y^i$ is a $1$-form on $M$ with $||\beta||_{\alpha}:=\sqrt{a^{ij}b_ib_j}<1$.
\end{definition}

	\begin{prof2}
		Suppose that $f_1$ is a nonzero constant function. Then, since $F_1$ and $F_2$ are Randers metrics, the third term of the right-hand side of the (\ref{6}) vanishes. Hence, the convolution metric becomes
		\begin{eqnarray}\label{35}
		_{f_2}F_1\star_{f_1}F_2
		=f_2^2(\alpha_1^2+2\alpha_1\beta_1+\beta_1^2)
		+f_1^2(\alpha_2^2+2\alpha_2\beta_2+\beta_2^2).
		\end{eqnarray}
		Since $f_1$ and $f_2$ are positive functions, denote $f_2\alpha_1, f_1\alpha_2,f_2\beta_1$ and $f_1\beta_2$ by $\alpha_1^*,\alpha_2^*, \beta_1^*$ and $\beta_2^*$ respectively. We obtain
		\begin{eqnarray}\label{36}
		_{f_2}F_1\star_{f_1}F_2
		=\alpha_1^{*2}+\alpha_2^{*2}+2(\alpha_1^{*}\beta_1^{*}+\alpha_2^{*}\beta_2^{*})
		+\beta_1^{*2}+\beta_2^{*2}.
		\end{eqnarray}
		By setting $\alpha=\sqrt{\alpha_1^{*2}+\alpha_2^{*2}}$ and 
		$\beta=\sqrt{\beta_1^{*2}+\beta_2^{*2}}$, by using the Definition \ref{defi33} and the relation (\ref{36}),  we get $\alpha\beta=\alpha_1^{*}\beta_1^{*}+\alpha_2^{*}\beta_2^{*}$. Hence the relation (\ref{34}) follows from the fact that $\alpha_1^{*}=f_2\alpha_1,\beta_1^{*}=f_2\beta_1, \alpha_2^{*}=f_1\alpha_2$ and $\beta_2^{*}=f_1\beta_2$. 
	\end{prof2}
\begin{corollary}
	A convolution Finslerian metric $_{f_2}F_1\star_{f_1}F_2$ is Euclidean if and only if $F_1$ as well as $F_2$ are Euclidean metrics and the functions $f_k,k\in\{1,2\}$ are both nonzero constants.
\end{corollary}

\section{Some examples of Finslerian convolution metrics}\label{Section4}
\begin{example}(Finslerian convolution metric of type Riemannian)
	Denote by $\| .   \|$ and $\langle ., .\rangle$ the standard Euclidean norm
	and the inner product in $\mathbb{R}^n$ respectively. Let $\mathbb{B}_1=\{x_1\in\mathbb{R}^3:\|x_1\|<1\}$ and $\mathbb{B}_2=\{x_2\in\mathbb{R}^2:\|x_2\|<1\}$ be two Euclidean balls. Define the functions 
	\begin{eqnarray}
	F_k:T\mathbb{B}_k&\longrightarrow&(0,\infty)\nonumber\\
	(x_k,y_k)&\to&\frac{\sqrt{\|y_k\|^2-(\|x_k\|^2\|y_k\|^2-\langle x_k,y_k\rangle^2)}}{1-\|x_k\|^2}
	\end{eqnarray}
	for $k\in\{1,2\}$. Then the function $F_1$ and $F_2$ are Riemannian metrics on $\mathbb{B}_1$ and $\mathbb{B}_2$ respectively known as Klein metrics.
	
	Consider two smooth functions $\rho_k:\mathbb{B}_k\longrightarrow \mathbb{R}$. Then the function 
	{\footnotesize{\begin{eqnarray}\label{4.2}
	F_{f_1,f_2}(x_1,x_2,y_1,y_2)
	&=&\left\{\frac{e^{2\rho_2}}{\left(1-\sum_{r=1}^3\big(x_1^r\big)^2\right)^2}
	\left[\sum_{r=1}^3\big(y_1^r\big)^2\left(1-\sum_{r=1}^3\big(x_1^r\big)^2\right)+\left(\sum_{r=1}^3x_1^ry_1^r\right)^2\right]\right.\nonumber\\
	&&+\frac{e^{2\rho_1}}{\left(1-\sum_{s=4}^5\big(x_2^s\big)^2\right)^2}
	\left[\sum_{s=4}^5\big(y_2^s\big)^2\left(1-\sum_{s=4}^5\big(x_2^s\big)^2\right)+\left(\sum_{s=4}^5x_2^sy_2^s\right)^2\right]\\
	&&+\left.2e^{2\rho_1}e^{2\rho_2}\sum_{r=1}^3\sum_{s=4}^5
	\frac{\partial \rho_1}{\partial x_1^r}
	\frac{\partial \rho_2}{\partial x_2^s}
	y_1^ry_2^s\right\}^{\frac{1}{2}}.\nonumber
	\end{eqnarray}}}
	is the Riemannian convolution of $F_1$ and $F_2$, via $f_1:=e^{\rho_1}$ and $f_2:=e^{\rho_2}$.
	
	When one the $\rho_k$ is a constant function or when $\nabla \rho_k \perp v$ for each tangent vector $v$ on $\mathbb{B}_k$, the Riemannian convolution function (\ref{4.2}) becomes a convolution Riemannian metric on $\mathbb{B}_1\times\mathbb{B}_2$.
\end{example}
\begin{example}(Finslerian convolution metric of type Minkowskian)
 Let consider the Finslerian function $F$ defined by the relation (\ref{dhggdgdgh}) in the Example \ref{example11}. 

It is easy to show that $F$ has the form $_{f_2}F_1\star_{f_1}F_2$, with $f_k:U\subset\mathbb{R}^2\longrightarrow (0,\infty)$ such that $f_k(x_1,x_2)=c_k$. That is $F$ is a Minkowskian convolution function on $M_1\times M_2$ where $M_1=\{(y^1,y^2)\in\mathbb{R}^2\}$, $M_2=\{(y^3,y^4)\in\mathbb{R}^2\}$, $F_1=\left(y_1^4+\lambda y_1^2y_2^2+y_2^4\right)^{\frac{1}{4}}$ and $F_2=
\left[y_3^2+y_4^2+\lambda \left(y_3^{2k}+y_4^{2k}\right)^{\frac{1}{k}}\right]^{\frac{1}{2}}$ via $c_1$ and $c_2$.
\end{example}
\begin{example}(Finslerian convolution metric of type Randers)
 Denote 
 by $\| .   \|_V$ the norm on a $3$-vector space $V$, by $\| .   \|$ the standard Euclidean norm and by $\langle ., .\rangle$ the inner product in $\mathbb{R}^{n-3}$. Let $\mathbb{B}=\{x_2\in\mathbb{R}^{n-3}:\|x_2\|<1\}$ be an Euclidean ball. 
	Then, for any $\epsilon\in [0,1)$, the following function $F:T\big(V\times \mathbb{B}\big)\longrightarrow (0,\infty)$ defined by 
 \begin{eqnarray}\label{dkkjfdjk}
    F(x_1,x_2,y_1,y_2)&=&F(x_1^1,x_1^2,x_1^3,x_2^4,...,x_2^n, y_1^1,y_1^2,y_1^3,y_2^4,...,y_2^n) \nonumber\\
    &=&\left\{\|x_2\|^2\left(\|y_1\|_V^2+\epsilon y_1^3\right)^2+2\sum_{\alpha=1}^3\sum_{\nu=4}^nx_1^{\alpha}x_2^{\nu}y_1^{\alpha}y_2^{\nu}\right.\nonumber\\
    &&+\left. \|x_1\|_V^2 \left[
    \frac{\sqrt{\|y_2\|^2-(\|x_2\|^2\|y_2\|^2-\langle x_2,y_2\rangle^2)}+ \langle x_2,y_2\rangle}{1-\|x_2\|^2}\right]^2
    \right\}^{\frac{1}{2}}
\end{eqnarray}
for any $x_1=(x_1^1,x_1^2,x_1^3)\in V$, $x_2=(x_2^4,...,x_2^n)\in\mathbb{B}, y_1=(y_1^1,y_1^2,y_1^3)\in T_{x_1}V\cong V$ and .$y_2=(y_2^4,...,y_2^n)\in T_{x_2}\mathbb{B}\cong\mathbb{R}^{n-3}$ such that 
{\small{$$
\|x_2\|^2\left(\|y_1\|_V^2+\epsilon y_1^3\right)^2
+\|x_1\|_V^2 \left[
    \frac{\sqrt{\|y_2\|^2-(\|x_2\|^2\|y_2\|^2-\langle x_2,y_2\rangle^2)}+ \langle x_2,y_2\rangle}{1-\|x_2\|^2}\right]^2
    >-2\sum_{\alpha=1}^3\sum_{\nu=4}^nx_1^{\alpha}x_2^{\nu}y_1^{\alpha}y_2^{\nu}.
$$}}
Then the function $F$ defined in (\ref{dkkjfdjk}) is Randers  convolution metric of $\|y_1\|_V^2+\epsilon y_1^3$ and $
    \frac{\sqrt{\|y_2\|^2-(\|x_2\|^2\|y_2\|^2-\langle x_2,y_2\rangle^2)}+ \langle x_2,y_2\rangle}{1-\|x_2\|^2}$, via $\|x_1\|$ and $\|x_2\|$.
\end{example}

\section*{Conflicts of Interest}
The author declares no conflicts of interest regarding the publication of this work.


\begin{thebibliography}{10}
\mark{RReferences}
 

\bibitem{bao}{D. Bao, S.-S. Chern and Z. Shen}, \textit{An Introduction to Riemann-Finsler Geometry}, 
Springer-Verlang New York, \textbf{200}, (2000), 1-417.
     
     \bibitem{bang2002} B.-Y. Chen, {\em Convolution of Riemannian manifolds and its applications}, 
Bull. Austral. Math. Soc. \textbf{66} (2002) 177-191.
%

\bibitem{bang2016}
   B.-Y. Chen, {\em Differential Geometry of
Warped Product Manifolds
and Submanifolds}, 
 	World Scientific \textit{ISBN 978-981-3208-92-6}, (2016) 1-78.
%

 \bibitem{Kuhnel}
 W. K\"uhnel and H.-B. Rademacher, {\em Conformal vector fields on pseudo-Riemannian spaces}, Diff. Geom. Appl. \textbf{7} (1997) 237-250.

	\bibitem{Nibaruta9} G. Nibaruta, M. Karimumuryango, A. Nibirantiza and D. Ndayirukiye, {\em 
	Twisted products Berwaldian metrics of polar type}, Differ. Geom.-Dyn. Syst. \textbf{22} (2020), to apper.
	%
 \bibitem{Bin}
 Y.-B Shen and Z. Shen, {\em Introduction to
Modern Finsler Geometry}, 
 Higher Education Press Limited Company and
World Scientific Publishing Co. Pte. Ltd (2016) 3-38.
\end{thebibliography}
 \end{document}